\documentclass[a4paper,10pt]{article}
\usepackage[latin1]{inputenc}
\usepackage[T1]{fontenc}
\usepackage{amsmath}
\usepackage{amsfonts}
\usepackage{amstext}
\usepackage{amsthm}
\usepackage[all]{xy}
\usepackage{geometry}
\usepackage{srcltx}
\usepackage{graphics}
\usepackage{epsfig}
\usepackage{subfigure}
\usepackage{slashed}
\usepackage{color}
\usepackage{amssymb}
\usepackage{srcltx}
\newtheorem{theorem}{Theorem}[section]

\DeclareMathOperator{\sfl}{sf}
\DeclareMathOperator{\diag}{diag}

\DeclareMathOperator{\sgn}{sgn}

\DeclareMathOperator{\Mat}{Mat}

\title{A Comparison Principle for Bifurcation of Periodic Solutions of Hamiltonian Systems}
 
\author{Helene Cyris, Joanna Janczewska and Nils Waterstraat}

\begin{document}
\date{}
\maketitle


\begin{abstract}
We obtain novel criteria for the existence of local bifurcation for periodic solutions of Hamiltonian systems by a comparison principle of the spectral flow. Our method allows to find the appearance of new solutions by a simple inspection of the coefficients of the system.  
\end{abstract}

\section{Introduction}
Let $\mathcal{H}:I\times\mathbb{R}\times\mathbb{R}^{2n}\rightarrow\mathbb{R}$ be a $C^2$-map, where $I=[0,1]$ denotes the unit interval, and assume that $\mathcal{H}(\lambda,t,u)$ is $2\pi$-periodic with respect to $t$ as well as $\nabla_u\mathcal{H}(\lambda,t,0)=0$ for all $(\lambda,t)\in I\times\mathbb{R}$. We consider the Hamiltonian systems

\begin{equation}\label{equation}
\left\{
\begin{aligned}
J u'(t)+\nabla_u&\mathcal{H}(\lambda,t,u(t))=0,\quad t\in [0,2\pi]\\
u(0)&=u(2\pi),
\end{aligned}
\right.
\end{equation}  
where $J$ denotes the standard symplectic matrix 

\begin{align}\label{J}
J=\begin{pmatrix}
0&-I_n\\
I_n&0
\end{pmatrix},
\end{align}
and note that $u\equiv 0$ is a solution of \eqref{equation} for all $\lambda\in I$ under the given assumptions. We say that $\lambda_0\in I$ is a bifurcation point for \eqref{equation}, if in every neighbourhood of $(\lambda_0,0)$ in $I\times C^1(S^1,\mathbb{R}^{2n})$ there is an element $(\lambda,u)$ such that $u\neq 0$ is a solution of \eqref{equation}. There have been numerous attempts to study the existence of bifurcation points in this setting (cf., e.g., \cite{Rabinowitz}, \cite{Mawhin}, \cite{Bartsch} and the references therein), and here we focus on an approach via the spectral flow that was introduced in \cite{Specflow} and applied to \eqref{equation} in \cite{SFLPejsachowiczII} and \cite{BifJac}.\\
The spectral flow is defined for any path $L=\{L_\lambda\}_{\lambda\in I}$ of selfadjoint Fredholm operators on a Hilbert space and, roughly speaking, it counts the net number of negative eigenvalues of $L_0$ that become positive whilst the parameter $\lambda$ traverses the interval $I$. In particular, a non-vanishing spectral flow means that there is an eigenvalue of $L_0$ that crosses $0$ for some $\lambda$ in $I$. The application of the spectral flow in bifurcation theory can be outlined as follows. Let $H$ be a real Hilbert space and $f:I\times H\rightarrow\mathbb{R}$ a $C^2$-map such that $0\in H$ is a critical point of all $f_\lambda:H\rightarrow\mathbb{R}$, i.e., $\nabla f_\lambda(0)=0$, $\lambda\in I$. We say that $\lambda^\ast\in I$ is a bifurcation point of the equation

\begin{align}\label{bifequ}
\nabla f_\lambda(u)=0
\end{align}
if in any neighbourhood of $(\lambda^\ast,0)$ in $I\times H$ there is some solution $(\lambda,u)$ of \eqref{bifequ} such that $u\neq 0$. The Hessian $L_\lambda$ of $f_\lambda$ at the critical point $0$ is a selfadjoint operator on $H$, and $L_{\lambda^\ast}$ has to be non-invertible if $\lambda^\ast$ is a bifurcation point of \eqref{bifequ}. If we require in addition that the operators $L_\lambda$ are Fredholm, then the spectral flow of the path $L=\{L_\lambda\}_{\lambda\in I}$ is defined and its non-triviality indicates that there is some $\lambda^\ast\in I$ such that $L_{\lambda^\ast}$ has a non-trivial kernel made by an eigenvalue that crosses $0$ when $\lambda$ passes $\lambda^\ast$. The main theorem of \cite{Specflow} is as follows:

\begin{theorem}\label{BifJacTheo}
Let $H$ be a real separable Hilbert space and $f:I\times H\rightarrow\mathbb{R}$ a $C^2$-map such that $0\in H$ is a critical point of all functionals $f_\lambda:H\rightarrow\mathbb{R}$. Let the Hessians $L_\lambda$ of $f_\lambda$ at $0$ be Fredholm operators and assume that $L_0$ and $L_1$ are invertible. If $\sfl(L)\neq 0\in\mathbb{Z}$ for the path $L=\{L_\lambda\}_{\lambda\in I}$, then there is a bifurcation point for \eqref{bifequ} in $(0,1)$. 
\end{theorem} 
Let us note that the spectral flow of $L$ is just the difference of the Morse indices of $L_0$ and $L_1$ if the latter are finite. Thus Theorem \ref{BifJacTheo} yields the well-known existence of a bifurcation point when the Morse indices are finite and change from $L_0$ to $L_1$ (cf. \cite{Conley}, \cite{Mawhin}).\\ 
Since the work \cite{Specflow} the spectral flow has been applied to study bifurcation problems for various kinds of differential equations with infinite Morse indices by several authors. Unfortunately, the spectral flow often is notoriously difficult to compute if possible at all. Pejsachowicz and the third author invented in \cite{BifJac} a comparison principle that can be used to estimate the spectral flow, which can be good enough to show that it does not vanish and so makes Theorem \ref{BifJacTheo} applicable (cf. \cite{Edinburgh}). For example, it allows to find bifurcation points for \eqref{equation} just from the eigenvalues of the family of Hessians of the Hamiltonian $\mathcal{H}$ at $0\in\mathbb{R}^{2n}$, which we now want to outline.\\
Let $H^\frac{1}{2}(S^1,\mathbb{R}^{2n})$ be the linear space of all functions $u:[0,2\pi]\rightarrow\mathbb{R}^{2n}$ such that
\begin{align}\label{series}
u(t)=c_0+\sum^\infty_{k=1}{(a_k\sin(kt)+b_k\cos(kt))},
\end{align}
where $c_0,a_k,b_k\in\mathbb{R}^{2n}$, $k\in\mathbb{N}$, and 
%\begin{align*}
$\sum^\infty_{k=1}{k(|a_k|^2+|b_k|^2)}<\infty$ (cf. \cite{AlbertoBuch}).\\
%\end{align*}
This is a Hilbert space with respect to the scalar product
\begin{align}\label{scalprod}
\langle u,v\rangle_{H^\frac{1}{2}}=2\pi \langle c_0,\tilde{c}_0\rangle+\pi\sum^\infty_{k=1}{k(\langle a_k,\tilde{a}_k\rangle+\langle b_k,\tilde{b}_k\rangle)},
\end{align}
where $\tilde{c}_0$ and $\tilde{a}_k,\tilde{b}_k$ are the Fourier coefficients of $v\in H^\frac{1}{2}(S^1,\mathbb{R}^{2n})$ as in \eqref{series}. \\
Let now
\begin{align*}
Q:H^\frac{1}{2}(S^1,\mathbb{R}^{2n})\times H^\frac{1}{2}(S^1,\mathbb{R}^{2n})\rightarrow\mathbb{R}
\end{align*}
be the continuous extension of the bounded bilinear form

\begin{align}\label{Qdiff}
\tilde{Q}(u,v)=\int^{2\pi}_0{\langle J u'(t),v(t)\rangle\,dt}
\end{align}
from the dense subspace $H^1(S^1,\mathbb{R}^{2n})$ to $H^\frac{1}{2}(S^1,\mathbb{R}^{2n})$, where we recall that $H^1(S^1,\mathbb{R}^{2n})$ consists of all absolutely continuous functions $u:[0,2\pi]\rightarrow\mathbb{R}^{2n}$ having a square integrable derivative. We consider the map

\begin{align}\label{f}
f:I\times H^\frac{1}{2}(S^1,\mathbb{R}^{2n})\rightarrow\mathbb{R},\quad f_\lambda(u)=\frac{1}{2}\,Q(u,u)+\int^{2\pi}_0{\mathcal{H}(\lambda,t,u(t))\,dt}.
\end{align} 
It is a well-known result (see \cite{Rabinowitz}, \cite{book}) that $f$ is $C^2$ if there are constants $a,b\geq 0$ and $r>1$ such that
	 \begin{align*}
\begin{split}
|\nabla_u\mathcal{H}(\lambda,t,u)|&\leq a+b|u|^r,\\
|D_u\nabla_u\mathcal{H}(\lambda,t,u)|&\leq a+b|u|^r,\quad (\lambda,t,u)\in  I\times\mathbb{R}\times\mathbb{R}^{2n}.
\end{split}
\end{align*}
Moreover, 
\[\langle \nabla_u f_\lambda ,v\rangle_{H^\frac{1}{2}}=Q(u,v)+\int^{2\pi}_0{\langle\nabla_u\mathcal{H}(\lambda,t,u(t)),v(t)\rangle\,dt},\quad v\in H^\frac{1}{2}(S^1,\mathbb{R}^{2n}),\]
and thus the critical points of $f_\lambda$ are the weak solutions of the Hamiltonian system \eqref{equation}. In particular,$ \ u\equiv 0\in H^\frac{1}{2}(S^1,\mathbb{R}^{2n})$ is a critical point of  $f_\lambda,$ and the associated Hessian is given by

\begin{align}\label{Hessian}
\langle L_\lambda u,v\rangle_{H^\frac{1}{2}}=Q(u,v)+\int^{2\pi}_0{\langle A_\lambda(t) u(t),v(t)\rangle\,dt},
\end{align}
where $A_\lambda(t)$ is the Hessian matrix of the functional $\mathcal{H}_\lambda(t,\cdot):\mathbb{R}^{2n}\rightarrow\mathbb{R}$ at $0$. Note that the kernel of $L_\lambda$ consists of the solutions of the linear equation

\begin{equation}\label{equationlin}
\left\{
\begin{aligned}
J u'(t)&+A_\lambda(t)u(t)=0,\quad t\in [0,2\pi]\\
u(0)&=u(2\pi).
\end{aligned}
\right.
\end{equation}  
The Hessians $L_\lambda$ are Fredholm operators (cf. \cite[Lemma 3.1]{CalcVar}, \cite{Fredholm}) and thus we see from Theorem \ref{BifJacTheo} that there is a bifurcation point of \eqref{equation} if \eqref{equationlin} only has the trivial solution for $\lambda=0,1$ and $\sfl(L)\neq 0$. As said before, obtaining the number $\sfl(L)$ usually is a challenging task (cf. \cite{SFLPejsachowiczII}) but estimating its value can be enough to show its non-triviality and thus to make Theorem \ref{BifJacTheo} applicable. This was done in a joint work of Pejsachowicz and the third author \cite{BifJac} as follows (cf. \cite{JSW}).\\  
We let $\mu_1(\lambda,t)\leq\cdots\leq\mu_{2n}(\lambda,t)$ be the eigenvalues of the matrices $A_\lambda(t)$ for $\lambda\in I$, $t\in[0,2\pi]$, and set

\begin{align}\label{alpha}
\alpha_\lambda:=\inf_{t\in[0,2\pi]}{\mu_1(\lambda,t)}
\end{align}
as well as

\begin{align}\label{beta}
\beta_\lambda:=\sup_{t\in[0,2\pi]}{\mu_{2n}(\lambda,t)}.
\end{align}
Thus $\alpha_\lambda$ is a lower bound for the spectra of the matrix family $\{A_\lambda(t)\}_{t\in[0,2\pi]}$ and $\beta_\lambda$ is an upper bound.

\begin{theorem}\label{thm-main}
Let $\mathcal{H}:I\times\mathbb{R}\times\mathbb{R}^{2n}\rightarrow\mathbb{R}$ be a family of Hamiltonians as above as well as $\alpha_1,\beta_0$ as in \eqref{alpha} and \eqref{beta}. Assume that \eqref{equationlin} only has the trivial solution for $\lambda=0$ and $\lambda=1$.

\begin{enumerate}
 \item[(i)] If there is some $k\in\mathbb{Z}$ such that  $\beta_0<k\leq\alpha_1$, then there is a bifurcation point of \eqref{equation}.
 \item[(ii)] If there are only finitely many $\lambda\in I$ such that \eqref{equationlin} has a non-trivial solution, then the number of bifurcation points of \eqref{equation} is greater or equal to the cardinality of 
 
 \[(\beta_0,\alpha_1]\cap\mathbb{Z}.\]
\end{enumerate}
\end{theorem}
Note that the assumption $\beta_0<k<\alpha_1$ is heavily restrictive as it requires the largest eigenvalue of the matrix family $A_0(t)$, $t\in[0,2\pi]$, to be strictly smaller than the smallest eigenvalue of $A_1(t)$, $t\in[0,2\pi]$. The aim of this note is to generalise Theorem \ref{thm-main} along the following idea. If we consider the diagonal matrices $C_0=\diag(\beta_0,\ldots,\beta_0)$ and $C_1=\diag(\alpha_1,\ldots,\alpha_1)$, then we in particular obtain from the assumption $\beta_0<k<\alpha_1$ the matrix inequality

\begin{align}\label{estimate}
A_0(t)\leq C_0\leq k I_{2n}\leq C_1\leq A_1(t),\quad\text{for all}\quad t\in[0,2\pi],
\end{align}
where $I_{2n}$ denotes the unit matrix and $\leq$ is the usual L\"owner-order on the set of all symmetric matrices (cf. \eqref{Loewner} below). The novelty of our work is to drop the middle term $k I_{2n}$ in \eqref{estimate} and to consider general symmetric matrices $C_0\leq C_1$ that are still bounded from below by $A_0(t)$ and from above by $A_1(t)$ for all $t\in[0,2\pi]$. In the subsequent second section we state two theorems that considerably strengthen Theorem \ref{thm-main}. The third section contains the proofs of these theorems, where we also provide a brief survey on the spectral flow and some methods to compute it.

%%%%%%%%%%%%%%%%%%%%%%%%%%%%%%%%%%%%%%%%%%%%%%%%%%%%%%%%%%%%%%%%%%%%%%%%%%%%%%%%%%%%%%%%%%%%%%%%%%%%%%%%%%%%%%%%%%%%%%%%%%%%%%%%%%%%%%%%%%%%%%%%%%%%%%%%%%%%%%%%%%%%%%%%%%%%%%%%%%%%%%%%%%%%%%%%%%%%%%%%%%%%%%%%%%%%%%%%%%%%%%%%%%%%%%%%%%%%%%%%%%%%%%%%%%%%%%%%%%%%%%%%%%%%%%%%%%%%%%%%%%%%%%%%%%%%%%%%%%%%%%%%%%%%%%%%%%%%%%%%%%%%%%%%%%%%%%%%%%%%%%%%%%%%%%%%%%%%%%%%%%%%%%%%%%%%%%%%%%%%%%%%%%%%%%%%%%%%%%%%%%%%%%%%%%%%%%%%%%%%%%%%%%%%%%%%%%%%%%%%%%%%

\section{The Main Theorems}
Henceforth we let $C_0,C_1\in\Mat(2n,\mathbb{R})$ be two symmetric matrices such that

\begin{align}\label{estimate1}
A_0(t)\leq C_0\leq C_1\leq A_1(t),\quad\text{for all}\quad t\in[0,2\pi]
\end{align}
and we denote their eigenvalues by

\[\mu_1(C_0)\leq\cdots\leq\mu_{2n}(C_0)\quad\text{as well as}\quad\mu_1(C_1)\leq\cdots\leq\mu_{2n}(C_1).\]
Note that geometrically, \eqref{estimate1} requires the paths $\{A_0(t)\}_{t\in[0,2\pi]}$ and $\{A_1(t)\}_{t\in[0,2\pi]}$ to run in affine convex cones with vertices made by $C_0$ or $C_1$, respectively.\\
Our first theorem on the existence of bifurcation points for \eqref{equation} is as follows.

\begin{theorem}\label{thm-mainI}
Let $\mathcal{H}:I\times\mathbb{R}\times\mathbb{R}^{2n}\rightarrow\mathbb{R}$ be a family of Hamiltonians as above and assume that \eqref{equationlin} only has the trivial solution for $\lambda=0$ and $\lambda=1$. Let $C_0,C_1\in\Mat(2n,\mathbb{R})$ be two symmetric matrices as in \eqref{estimate1}. If there is some $1\leq i\leq 2n$ such that 

\begin{align}\label{muequ}
\mu_i(C_0)<0\leq\mu_i(C_1),
\end{align}
then there is a bifurcation point of periodic solutions from the trivial branch for \eqref{equation}.
\end{theorem}
Note that \eqref{muequ} just means that the Morse index $\mu_-(C_0)$ of the matrix $C_0$, i.e. the number of negative eigenvalues including multiplicities, differs from the Morse index $\mu_-(C_1)$ of $C_1$.\\
Let us emphasize that the existence of a bifurcation point for \eqref{equation} follows from Theorem \ref{BifJacTheo} if $\sfl(L)\neq 0$, where the path $L=\{L_\lambda\}_{\lambda\in I}$ is made by the operators \eqref{Hessian} in the Hilbert space $H^\frac{1}{2}(S^1,\mathbb{R}^{2n})$. This spectral flow will in general be hard to compute if possible at all. Theorem \ref{thm-mainI} shows that if there are two matrices $C_0, C_1$ as in \eqref{estimate1}, then bifurcation follows from the non-triviality of the number $\mu_-(C_0)-\mu_-(C_1)$. The latter actually is the spectral flow of the path $I\ni\lambda\mapsto(1-\lambda)C_0+\lambda C_1\in\Mat(2n,\mathbb{R})$ of symmetric matrices in the finite dimensional Hilbert space $\mathbb{R}^{2n}$ (see \cite[\S 1]{book}). Thus to some extend Theorem \ref{thm-mainI} reduces the computation of the spectral flow $\sfl(L)$ to finite dimensions.\\
In a second theorem we wish to further strengthen Theorem \ref{thm-mainI} under the additional assumption that $C_0$ and $C_1$ commute with the symplectic standard matrix $J$ in \eqref{J}. Let us note that a symmetric matrix in $\Mat(2n,\mathbb{R})$ commutes with $J$ if and only if it is of the form

\begin{align*}
\begin{pmatrix}
A&B\\
-B&A
\end{pmatrix}
\end{align*}
for some $A=A^T\in\Mat(n,\mathbb{R})$ and $B=-B^T\in\Mat(n,\mathbb{R})$. Our second theorem is as follows.

\begin{theorem}\label{thm-mainII}
Let $\mathcal{H}:I\times\mathbb{R}\times\mathbb{R}^{2n}\rightarrow\mathbb{R}$ be a family of Hamiltonians as above and assume that \eqref{equationlin} only has the trivial solution for $\lambda=0$ and $\lambda=1$. Let $C_0,C_1\in\Mat(2n,\mathbb{R})$ be two symmetric matrices as in \eqref{estimate1} and assume in addition that $[J,C_0]=[J,C_1]=0$. 

\begin{itemize}
 \item[(i)] If there are some $i\in\{1,\ldots,2n\}$ and $k\in\mathbb{Z}$ such that 

\begin{align}\label{thmIIcond}
\mu_i(C_0)<k\leq\mu_i(C_1),
\end{align}
then there is a bifurcation point of periodic solutions from the trivial branch for \eqref{equation}.
  \item[(ii)] If there are only finitely many $\lambda\in I$ such that \eqref{equationlin} has a non-trivial solution, then the number of bifurcation points of \eqref{equation} is at least
  
  \[\frac{\sum^{2n}_{i=1}|(\mu_i(C_0),\mu_i(C_1)]\cap\mathbb{Z}|}{2n},\]
where $|\cdot|$ stands for the cardinality of a set.
\end{itemize}
\end{theorem}
Note that $\mu_i(C_0)\leq\mu_i(C_1)$, $1\leq i\leq 2n$, by \eqref{estimate1}, and \eqref{thmIIcond} means that there is a pair of eigenvalues that is separated by an integer. For $k=0$, Theorem \ref{thm-mainII} (i) is a direct consequence of Theorem \ref{thm-mainI}, but for $k\neq 0$ the additional assumption $[J,C_0]=[J,C_1]=0$ on the matrices $C_0$, $C_1$ in Theorem \ref{thm-mainII} will be essential in its proof. Note that the latter is in particular satisfied for the matrices $C_0=\diag(\beta_0,\ldots,\beta_0)$ and $C_1=\diag(\alpha_1,\ldots,\alpha_1)$ in \eqref{estimate} and thus Theorem \ref{thm-main} is an immediate consequence of Theorem \ref{thm-mainII}.

\section{Proof of Theorem \ref{thm-mainI} and Theorem \ref{thm-mainII}}
We now firstly recap some basics on the spectral flow and make a first common step in the proofs Theorem \ref{thm-mainI} and Theorem \ref{thm-mainII}. Afterwards we deal with them in two individual sections.

%%%%%%%%%%%%%%%%%%%%%%%%%%%%%%%%%%%%%%%%%%%%%%%%%%%%%%%%%%%%%%%%%%%%%%%%%%%%%%%%%%%%%%%%%%%%%%%%%%%%%%%%%%%%%%%%%%%%%%%%%%%%%%%%%%%%%%%%%%%%%%%%%%%%%%%%%%%%%%%%%%%%%%%%%%%%%%%%%%%%%%%%%%%%%%%%%%%%%%%%%%%%%%%%%%%%%%%%%%%%%%%%%%%%%%%%%%%%%%%%%%%%%%%%%%%%%%%%%%%%%%%%%%%%%%%%%%%%%%%%%%%%%%%%%%%%%%%%%%%%%%%%%%%%%%%%%%%%%%%%%%%%%%%%%%%%%%%%%%%%%%%%%%%%%%%%%%%%%%%%%%%%%%%%%%%%%%%%%%%%%%%%%%%%%%%%%%%%%%%%%%%%%%%%%%%%%%%%%%%%%%%%%%%%%%%%%%%%%%%%%%%%%%%%%%%%%%%%%%%%%%%%%%%%%%%%%%%%%%%%%%%%%%%%%%%%%%%%%%%%%%%%%%%%%%%%%%%%%%%%%%%%%%%%%%%%%%%%%%%%%%%%%%%%%%%%%%

\subsection{The Spectral Flow Setting}
Let $C_0,C_1\in\Mat(2n,\mathbb{R})$ be two symmetric matrices as in \eqref{estimate1}, i.e.

\begin{align*}
A_0(t)\leq C_0\leq C_1\leq A_1(t),\quad\text{for all}\quad t\in[0,2\pi],
\end{align*}
and set

\[C_\lambda:=(1-\lambda)C_0+\lambda C_1\in\Mat(2n,\mathbb{R}),\quad\lambda\in[0,1].\]
We consider the homotopy $H=\{H_{(\lambda,s)}\}_{(\lambda,s)\in I\times I}$ of selfadjoint operators that is uniquely determined by

\begin{align*}
\langle H_{(\lambda,s)} u,v\rangle_{H^\frac{1}{2}}=Q(u,v)+\int^{2\pi}_0{\langle ((1-s)A_\lambda(t)+s C_\lambda) u(t),v(t)\rangle\,dt},
\end{align*}
and note that each $H_{(\lambda,s)}$ is Fredholm by, e.g., \cite{Fredholm}, \cite[Lemma 3.1]{CalcVar}. Thus the spectral flows of the paths obtained by restricting the parameters of $H=\{H_{(\lambda,s)}\}_{(\lambda,s)\in I\times I}$ to any edge of the rectangle $I\times I$ are defined, and we obtain

\[\sfl(H_{(\cdot,0)})=\sfl(H_{(0,\cdot)})+\sfl(H_{(\cdot,1)})-\sfl(H_{(1,\cdot)})\]
as the spectral flow is homotopy invariant, it vanishes for constant paths, and changes its sign when traversing a path in opposite direction (cf., e.g., \cite{CompSfl}). Another important property of the spectral flow comes from the well-known partial order

\begin{align}\label{Loewner}
A\leq B:\Longleftrightarrow \langle (B-A)u,u\rangle\geq 0,\,\text{for all}\, u
\end{align}
on the set of selfadjoint operators on a Hilbert space. A path $\{L_\lambda\}_{\lambda\in I}$ of selfadjoint Fredholm operators is called increasing/decreasing if $L_\lambda\leq L_\mu$/ $L_\mu\leq L_\lambda$ whenever $\lambda\leq\mu$. It was shown in \cite{BifJac} (cf. \cite{JSW}) that $\sfl(L)\geq 0$ for increasing paths, and $\sfl(L)\leq 0$ for decreasing ones. As it follows from \eqref{estimate1} that the path made by

\begin{align*}
\langle H_{(0,s)} u,v\rangle_{H^\frac{1}{2}}&=Q(u,v)+\int^{2\pi}_0{\langle ((1-s)A_0(t)+s C_0) u(t),v(t)\rangle\,dt},\quad s\in I,\\
\end{align*}
is increasing, and the one made by the operators

\begin{align*}
\langle H_{(1,s)} u,v\rangle_{H^\frac{1}{2}}=Q(u,v)+\int^{2\pi}_0{\langle ((1-s)A_1(t)+s C_1) u(t),v(t)\rangle\,dt},\quad s\in I,
\end{align*}
is decreasing, we see that 

\begin{align}\label{estimateL}
\sfl(L)=\sfl(H(\cdot,0))=\sfl(H_{(0,\cdot)})+\sfl(H_{(\cdot,1)})-\sfl(H_{(1,\cdot)})\geq \sfl(H_{(\cdot,1)})\geq 0,
\end{align}
where the latter inequality follows once again by \eqref{estimate1}. Henceforth we abbreviate 

\[M_\lambda:=H_{(\lambda,1)},\quad \lambda\in I,\]
and note that

\begin{align*}
\langle M_\lambda u,v\rangle_{H^\frac{1}{2}}=Q(u,v)+\int^{2\pi}_0{\langle C_\lambda u(t),v(t)\rangle\,dt}.
\end{align*}
Hence by \eqref{estimateL}

\begin{align}\label{LgeqM}
\sfl(L)\geq\sfl(M)
\end{align}
and by Theorem \ref{BifJacTheo} there is a bifurcation point for \eqref{equation} if $\sfl(M)>0$, where $M=\{M_\lambda\}_{\lambda\in I}$. Let us note that Theorem \ref{BifJacTheo} only shows that there is a bifurcation of weak solutions of \eqref{equation} in $I\times H^\frac{1}{2}(S^1,\mathbb{R}^{2n})$. It is well-known that the critical points of the functionals $f_\lambda$ in \eqref{f} actually are classical solutions of \eqref{equation} which converge to $0$ in $C^1(S^1,\mathbb{R}^{2n})$ under the given assumptions, and detailed proofs can be found, e.g., in \cite{Rabinowitz} or \cite[Thm. 12.3.3]{book}.\\ 
In the following two sections we show that $\sfl(M)>0$ under the assumptions of either Theorem \ref{thm-mainI} or Theorem \ref{thm-mainII}, where we use a common method to compute the spectral flow that was invented by Robbin and Salamon in \cite{Robbin-Salamon}.  Let $\mathcal{L}=\{\mathcal{L}_\lambda\}_{\lambda\in [a,b]}$ be a continuously differentiable path in the normed space of bounded operators on a Hilbert space, where we assume in addition that each $\mathcal{L}_\lambda$ is selfadjoint and Fredholm. Then $\lambda^\ast\in[a,b]$ is called \textit{crossing} of $\mathcal{L}$ if $\ker(\mathcal{L}_{\lambda^\ast})\neq \{0\}$. The \textit{crossing form} of a crossing $\lambda^\ast$ is the quadratic form on the finite dimensional space $\ker(\mathcal{L}_{\lambda^\ast})$ given by

\[\Gamma(\mathcal{L},\lambda^\ast)[u]=\langle \dot{\mathcal{L}}_{\lambda^\ast}u,u\rangle,\qquad u\in \ker(\mathcal{L}_{\lambda^\ast}),\]   
where $\dot{\mathcal{L}}_{\lambda^\ast}$ denotes the derivative of the path $\mathcal{L}$ at $\lambda=\lambda^\ast$. Finally, a crossing is called \textit{regular} if $\Gamma(\mathcal{L},\lambda^\ast)$ is non-degenerate. The following theorem was shown in \cite{Robbin-Salamon} (see also \cite{Homoclinics}).

\begin{theorem}\label{thm-sfl-crossings}
Let $\mathcal{L}=\{\mathcal{L}_\lambda\}_{\lambda\in[a,b]}$ be a continuously differentiable path of selfadjoint Fredholm operators as above.
\begin{enumerate}
\item[(i)] If $\mathcal{L}$ has only regular crossings, then there are only finitely many of them and the spectral flow of $\mathcal{L}$ is given by

\begin{align}\label{sfl-formula}
\sfl(\mathcal{L},[a,b])=-m^-(\Gamma(\mathcal{L},a))+\sum_{\lambda\in(a,b)}\sgn(\Gamma(\mathcal{L},\lambda))+m^-(-\Gamma(\mathcal{L},b)),
\end{align} 
where $m^-$ denotes the Morse-index of a quadratic form and $\sgn$ its signature.
\item[(ii)] There is some $\varepsilon>0$ such that the perturbed path $\mathcal{L}^\delta=\{\mathcal{L}_\lambda+\delta I_H\}_{\lambda\in[a,b]}$ has only regular crossings for almost every $\delta\in(-\varepsilon,\varepsilon)$.
\end{enumerate}
\end{theorem}
Finally, below we need the following perturbation invariance of the spectral flow, which can be found, e.g., in \cite[Lemma 2.1]{IJW}.

\begin{theorem}\label{thm-sfl-perturbation}
For every path $\mathcal{L}=\{\mathcal{L}_\lambda\}_{\lambda\in[a,b]}$ of selfadjoint Fredholm operators there is some $\varepsilon>0$ such that

\[\sfl(\mathcal{L})=\sfl(\mathcal{L}^\delta),\]
for all $\delta\in[0,\varepsilon)$, where $\mathcal{L}^\delta=\{\mathcal{L}_\lambda+\delta I_H\}_{\lambda\in[a,b]}$. 
\end{theorem}

We now prove Theorem \ref{thm-mainI} and Theorem \ref{thm-mainII} in two separate sections. Let us briefly introduce the strategy of both proofs. By Theorem \ref{thm-sfl-crossings} and Theorem \ref{thm-sfl-perturbation}, there is some $\varepsilon>0$ such that for almost all $\delta\in[0,\varepsilon)$ the path $M^\delta$ has only regular crossings and $\sfl(M^\delta)=\sfl(M)$. Moreover, $\sfl(M^\delta)$ is given by \eqref{sfl-formula}, where

\[\Gamma(M^\delta,\lambda^\ast)[u]=\Gamma(M,\lambda^\ast)[u]=\langle (C_1-C_0)u,u\rangle_{H^\frac{1}{2}},\quad u\in\ker(M^\delta).\]  
As $C_1\geq C_0$ by assumption, it follows from \eqref{scalprod} that $\Gamma(M^\delta,\lambda^\ast)$ is positive semidefinite and moreover, as $\Gamma(M^\delta,\lambda^\ast)$ is non-degenerate by our choice of $\delta$, it is positive definite. Thus by \eqref{LgeqM} and \eqref{sfl-formula}

\begin{align}\label{sfl-formulaII}
\sfl(L)\geq\sfl(M)=\sfl(M^\delta)=\sum_{\lambda\in(0,1]}\dim\ker(M^\delta_\lambda),
\end{align} 
which finally shows that $\sfl(L)\neq0$ if we just can find some $\lambda\in(0,1]$ such that $\ker(M^\delta_\lambda)$ is non-trivial. Let us emphasize that \eqref{sfl-formulaII} holds for almost all $\delta\in[0,\varepsilon)$ and some $\varepsilon>0$.\\
We now consider Theorem \ref{thm-mainI} and Theorem \ref{thm-mainII} separately. In what follows we need in addition to the scalar product \eqref{scalprod} on $H^\frac{1}{2}(S^1,\mathbb{R}^{2n})$ also the scalar product

\begin{align}\label{scalprod2}
\langle u,v\rangle_{L^2}=2\pi \langle c_0,\tilde{c}_0\rangle+\pi\sum^\infty_{k=1}{\langle a_k,\tilde{a}_k\rangle+\langle b_k,\tilde{b}_k\rangle}
\end{align}
for $u,v\in L^2(S^1,\mathbb{R}^{2n})$ as in \eqref{series}.

%%%%%%%%%%%%%%%%%%%%%%%%%%%%%%%%%%%%%%%%%%%%%%%%%%%%%%%%%%%%%%%%%%%%%%%%%%%%%%%%%%%%%%%%%%%%%%%%%%%%%%%%%%%%%%%%%%%%%%%%%%%%%%%%%%%%%%%%%%%%%%%%%%%%%%%%%%%%%%%%%%%%%%%%%%%%%%%%%%%%%%%%%%%%%%%%%%%%%%%%%%%%%%%%%%%%%%%%%%%%%%%%%%%%%%%%%%%%%%%%%%%%%%%%%%%%%%%%%%%%%%%%%%%%%%%%%%%%%%%%%%%%%%%%%%%%%%%%%%%%%%%%%%%%%%%%%%%%%%%%%%%%%%%%%%%%%%%%%%%%%%%%%%%%%%%%%%%%%%%%%%%%%%%%%%%%%%%%%%%%%%%%%%%%%%%%%%%%%%%%%%%%%%%%%%%%%%%%%%%%%%%%%%%%%%%%%%%%%%%%%%%%

\subsection{Proof of Theorem \ref{thm-mainI}}
As $\mu_i(C_0)<0\leq\mu_i(C_1)$ by \eqref{muequ}, for every sufficiently small $\delta>0$, there is some $\lambda^\ast\in(0,1)$ such that $\mu_i(C_{\lambda^\ast})=-\delta$. Moreover, we can assume that \eqref{sfl-formulaII} holds for this $\delta$.\\
Let now $u_0\in\ker(C_{\lambda^\ast}+\delta I_{2n})$, where $I_{2n}\in\Mat(2n,\mathbb{R})$ is the identity matrix, and let $u\in H^\frac{1}{2}(S^1,\mathbb{R}^{2n})$ be the constant function $u(t)\equiv u_0$. Then

\begin{align*}
\langle M^\delta_{\lambda^\ast}u,v\rangle_{H^\frac{1}{2}}&=Q(u,v)+\delta \langle u,v\rangle_{H^\frac{1}{2}}+\int^{2\pi}_0{\langle C_{\lambda^\ast} u(t),v(t)\rangle\,dt}\\
&=\delta \langle u,v\rangle_{H^\frac{1}{2}}-\delta\int^{2\pi}_0{\langle u(t),v(t)\rangle\,dt}=0,\quad v\in H^\frac{1}{2}(S^1,\mathbb{R}^{2n}),
\end{align*}
where we use that $Q(u,v)=0$ by \eqref{Qdiff} and $\langle u,v\rangle_{H^\frac{1}{2}}=\langle u,v\rangle_{L^2}$ for every $v\in H^\frac{1}{2}(S^1,\mathbb{R}^{2n})$ by \eqref{scalprod} and \eqref{scalprod2} as $u$ is constant. Thus $u\in\ker M^\delta_{\lambda^\ast}$, which shows that $\sfl(L)\geq 1$ by \eqref{sfl-formulaII} and hence the existence of a bifurcation of periodic solutions for \eqref{equation} by Theorem \ref{BifJacTheo}. This completes the proof of Theorem \ref{thm-mainI}.

%%%%%%%%%%%%%%%%%%%%%%%%%%%%%%%%%%%%%%%%%%%%%%%%%%%%%%%%%%%%%%%%%%%%%%%%%%%%%%%%%%%%%%%%%%%%%%%%%%%%%%%%%%%%%%%%%%%%%%%%%%%%%%%%%%%%%%%%%%%%%%%%%%%%%%%%%%%%%%%%%%%%%%%%%%%%%%%%%%%%%%%%%%%%%%%%%%%%%%%%%%%%%%%%%%%%%%%%%%%%%%%%%%%%%%%%%%%%%%%%%%%%%%%%%%%%%%%%%%%%%%%%%%%%%%%%%%%%%%%%%%%%%%%%%%%%%%%%%%%%%%%%%%%%%%%%%%%%%%%%%%%%%%%%%%%%%%%%%%%%%%%%%%%%%%%%%%%%%%%%%%%%%%%%%%%%%%%%%%%%%%%%%%%%%%%%%%%%%%%%%%%%%%%%%%%%%%%%%%%%%%%%%%%%%%%%%%%%%%%%%%%%

\subsection{Proof of Theorem \ref{thm-mainII}}
As $(i)$ follows from Theorem \ref{thm-mainI} for $k=0$, we can henceforth assume that $k\neq 0$. We now begin similarly as in the proof of Theorem \ref{thm-mainI} and note that, as $\mu_i(C_0)<k\leq\mu_i(C_1)$ by \eqref{thmIIcond}, for every sufficiently small $\delta>0$, there is some $\lambda^\ast\in(0,1)$ such that $\mu_i(C_{\lambda^\ast})=k-\delta|k|$. Moreover, we can assume that \eqref{sfl-formulaII} holds for the number $\delta$.\\
Let $u_0\in\ker(C_{\lambda^\ast}-(k-\delta |k|) I_{2n})$. By assumption, $[J,C_0]=[J,C_1]=0$ and thus $[J,C_{\lambda^\ast}]=0$. Consequently, $Ju_0\in\ker(C_{\lambda^\ast}-(k-\delta |k|) I_{2n})$, and we now consider $u\in H^\frac{1}{2}(S^1,\mathbb{R}^{2n})$ given by

\begin{align}\label{repu}
u(t)=u_0 \cos(kt)+Ju_0 \sin(kt)=\sum^\infty_{j=1}{a_j \cos(jt)+b_j \sin(jt),}
\end{align}
where $a_{|k|}=u_0$, $b_{|k|}=\sgn(k) J u_0$, and $a_j=b_j=0$ for $j\neq |k|$. Then $C_{\lambda^\ast}u=(k-\delta |k|)\,u$, it follows from \eqref{scalprod} and \eqref{scalprod2} that

\[\langle u,v\rangle_{H^\frac{1}{2}}=|k|\langle u,v\rangle_{L^2},\quad v\in H^\frac{1}{2}(S^1,\mathbb{R}^{2n}),\]
and we obtain from \eqref{Qdiff} that

\[Q(u,v)=-k\langle u,v\rangle_{L^2},\quad v\in H^\frac{1}{2}(S^1,\mathbb{R}^{2n}).\]
Consequently,

\begin{align*}
\langle M^\delta_{\lambda^\ast}u,v\rangle&=Q(u,v)+\delta\langle u,v\rangle_{H^\frac{1}{2}}+\langle C_{\lambda^\ast} u,v\rangle_{L^2}\\
&=-k\langle u,v\rangle_{L^2}+\delta|k| \langle u,v\rangle_{L^2}+\langle (k-\delta |k|)u,v\rangle_{L^2}\\
&=0
\end{align*}
for any $v\in H^\frac{1}{2}(S^1,\mathbb{R}^{2n})$.  Thus $u\in\ker M^\delta_{\lambda^\ast}$, which again shows that $\sfl(L)\geq 1$ by \eqref{sfl-formulaII} and hence the existence of a bifurcation of periodic solutions for \eqref{equation} by Theorem \ref{BifJacTheo}, as claimed in part $(i)$ of Theorem \ref{thm-mainII}.\\ 
For part $(ii)$, we first recall the following result from \cite[Th. 2.1 (ii)]{BifJac}, where we use the notation from Theorem \ref{BifJacTheo}.

\begin{theorem}
If there are only finitely many $\lambda\in (0,1)$ such that $L_\lambda$ is non-invertible, then the family $f$ has at least

\[\frac{|\sfl(L)|}{\max\{\dim\ker L_\lambda:\, \lambda\in(0,1)\}}\]
bifurcation points in $(0,1)$.
\end{theorem}
As $\ker L_\lambda$ is made by the solutions of \eqref{equationlin}, we see that

\[\max\{\dim\ker L_\lambda:\, \lambda\in(0,1)\}\leq 2n,\]
and thus it remains to show that

\begin{align}\label{finalestimate}
\sfl(L)\geq\sum^{2n}_{i=1}|(\mu_i(C_0),\mu_i(C_1)]\cap\mathbb{Z}|.
\end{align}  
Let $\delta>0$ be such that \eqref{sfl-formulaII} holds for $M^\delta$. Moreover, we can assume that whenever $k\in(\mu_i(C_0),\mu_i(C_1)]\cap\mathbb{Z}$ for some $1\leq i\leq 2n$, then $k-\delta|k|\in(\mu_i(C_0),\mu_i(C_1)]$.\\ 
We let $k_1,\ldots, k_m\in\mathbb{Z}$ denote the elements of

\[(\mu_1(C_0),\mu_1(C_1)]\cap\mathbb{Z}.\]
As in the proof of the first part of the theorem, for each $k_i$ there is some $\lambda^1_i\in(0,1)$ such that $\mu_1(C_{\lambda^1_i})=k_i-\delta |k_i|$ and 

\[\dim\ker (M^\delta_{\lambda^1_i})\geq 1,\quad 1\leq i\leq m.\]
It follows from \eqref{sfl-formulaII} that

\[\sfl(L)\geq|(\mu_1(C_0),\mu_1(C_1)]\cap\mathbb{Z}|.\]
Let now $l_1,\ldots, l_n\in\mathbb{Z}$ denote the elements of 

\[(\mu_2(C_0),\mu_2(C_1)]\cap\mathbb{Z}.\]
Again for each $l_j$, there is some $\lambda^2_j\in(0,1)$ such that $\mu_2(C_{\lambda^2_j})=l_j-\delta |l_j|$. If $\lambda^1_i=\lambda^2_j$ for some $i,j$, then there are two cases. Firstly, if $\mu_1(C_{\lambda^1_i})=\mu_2(C_{\lambda^1_i})$, then the eigenspace of $C_{\lambda^1_i}$ for the eigenvalue $\mu_1(C_{\lambda^1_i})$ is at least of dimension two. Secondly, if $\mu_1(C_{\lambda^1_i})<\mu_2(C_{\lambda^1_i})$, then the corresponding eigenspaces intersect trivially. In both cases we obtain two linearly independent eigenvectors of $C_{\lambda^1_i}$ that yield linearly independent functions in \eqref{repu}. Consequently, $\dim\ker (M^\delta_{\lambda^1_i})\geq 2$, where one dimension was found by considering the elements of $(\mu_1(C_0),\mu_1(C_1)]\cap\mathbb{Z}$ and one dimension from $(\mu_2(C_0),\mu_2(C_1)]\cap\mathbb{Z}$. \\
In summary, after these first two steps of the construction, we see from \eqref{sfl-formulaII} that

\[\sfl(L)\geq|(\mu_1(C_0),\mu_1(C_1)]\cap\mathbb{Z}|+|(\mu_2(C_0),\mu_2(C_1)]\cap\mathbb{Z}|.\]   
If we continue this process until $(\mu_{2n}(C_0),\mu_{2n}(C_1)]\cap\mathbb{Z}$, we finally arrive at \eqref{finalestimate} which shows the second part of Theorem \ref{thm-mainII}.

\vspace*{1.3cm}

\begin{minipage}{1.2\textwidth}
\begin{minipage}{0.4\textwidth}

Helene Cyris\\
Martin-Luther-Universit\"at Halle-Wittenberg\\
Naturwissenschaftliche Fakult\"at II\\
Institut f\"ur Mathematik\\
06099 Halle (Saale)\\
Germany\\
helene.cyris@mathematik.uni-halle.de\\\\\\
Joanna Janczewska\\
Institute of Applied Mathematics\\
Faculty of Applied Physics and Mathematics\\
Gda\'{n}sk University of Technology\\
Narutowicza 11/12, 80-233 Gda\'{n}sk, Poland\\
joanna.janczewska@pg.edu.pl\\\\

\end{minipage}
\hfill
\begin{minipage}{0.6\textwidth}

Nils Waterstraat\\
Martin-Luther-Universit\"at Halle-Wittenberg\\
Naturwissenschaftliche Fakult\"at II\\
Institut f\"ur Mathematik\\
06099 Halle (Saale)\\
Germany\\
nils.waterstraat@mathematik.uni-halle.de
\end{minipage}
\end{minipage}

\end{document}